\documentclass[10pt,a4paper,oneside,onecolumn]{article}

\usepackage{apacite} % needs to be first, otherwise weird effects
\usepackage{ICROMA} % template by RailDresden
\usepackage{enumerate}
\usepackage{xspace}
\usepackage{accents}
\usepackage{amsfonts}
\usepackage{algorithm}
\usepackage[noend]{algpseudocode}

\title{Real-time Optimization of Transport Chains for Single Wagon Load Railway Transport}
\author{
    Carsten Moldenhauer$^1$,
    Philipp Germann$^1$,
    Cedric Heimhofer$^2$, \\
    Caroline Spieckermann, % no affiliation is intended and OK
    Andreas Andresen$^1$
}
\affiliation{
    $^1$Swiss National Railways SBB Cargo AG, Olten, Switzerland \\
   $^2$ Accenture AG, Zurich, Switzerland \\
   E-mail: \href{mailto:carsten.moldenhauer@sbbcargo.com}{carsten.moldenhauer@sbbcargo.com},
   \href{mailto:philipp.germann@sbbcargo.com}{philipp.germann@sbbcargo.com}
}

\newcommand{\WLV}{WLV\xspace}
\newcommand{\SBBCargo}{SBB Cargo\xspace}
\newcommand{\tc}{transport chain\xspace}
\newcommand{\tcs}{transport chains\xspace}
\newcommand{\eg}{e.g.\xspace}
\newcommand{\ie}{i.e.\xspace}
\newcommand{\sref}[1]{\hyperref[#1]{Section~\ref*{#1}}}
\newcommand{\fref}[1]{\hyperref[#1]{Figure~\ref*{#1}}}
\newcommand{\aref}[1]{\hyperref[#1]{Algorithm~\ref*{#1}}}
\newcommand{\TKModul}{TK-Modul\xspace}
\newcommand{\ubar}[1]{\underaccent{\bar}{#1}}

\makeatletter
\newcommand{\customlabel}[2]{%
   \protected@write \@auxout {}{\string \newlabel {#1}{{#2}{\thepage}{#2}{#1}{}} }%
   \hypertarget{#1}{#2}
}
\makeatother

\newlength{\bibparskip}
\let\oldthebibliography\thebibliography
\renewcommand\thebibliography[1]{%
  \oldthebibliography{#1}%
  \setlength{\parskip}{\bibitemsep}%
  \setlength{\itemsep}{\bibparskip}%
}

\newenvironment{varalgorithm}[1]
  {\algorithm\renewcommand{\thealgorithm}{#1}}
  {\endalgorithm}

\begin{document}
% Pagenumbering is removed (for proceeding generation)

\maketitle

% The abstract should not exceed 250 words.
\begin{abstract}
The freight branch of the Swiss national railways, \SBBCargo, offers customers to ship single or few wagons within its wagon load transportation system (\WLV). In this system, wagons travel along a \tc which is a sequence of consecutive trains.
Recently, \SBBCargo redesigned its IT systems and renewed the computation of these \tcs.
This paper describes the main design decisions and technical details: data structures, search algorithms, mathematical optimization of throughput in the real-time setting, and some selected details for making the algorithms work in the operational software. We also comment on the employed technology stack and finally demonstrate some performance metrics from running operations.
\end{abstract}

% Up to five keywords.

\section{Introduction}
\begin{figure}[t]
    \begin{minipage}[c]{0.55\textwidth}
        \includegraphics[width=\textwidth]{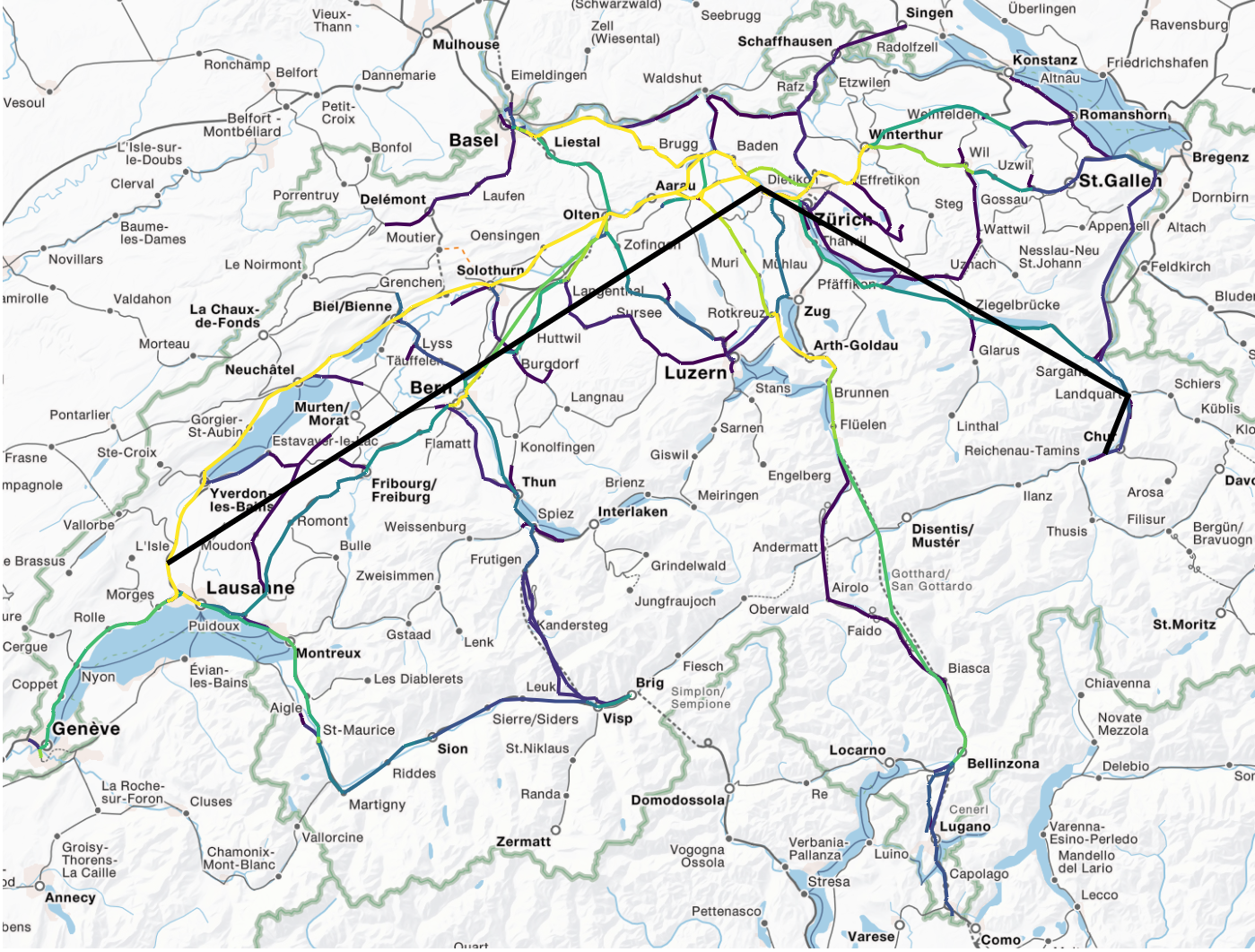}
    \end{minipage}\hfill
    \begin{minipage}[c]{0.4\textwidth}
        black: a transport chain from Felsberg via Landquart and Zurich to Cossonay.\\
        \\
        thick colored lines: \WLV network colored by average number of transported wagons in 2023 (blue is low, yellow is highest).\\
        \\
        thin grey lines: the rail network in Switzerland.
    \caption{The \WLV network of Switzerland with examplary \tc from Felsberg to Cossonay.}
    \end{minipage}
    \label{fig:transport_chain}
\end{figure}

Switzerland traditionally exhibits a strong modal split with a share of 35 to 40\% of all goods being transported on rails~\citep{FSO}. Of these railway transports, the freight branch of the swiss national railways \emph{\SBBCargo}, holds a market share of about 34\% which accounts for 26M net tons transported, or 4,622M net ton kilometers, in 2023~\citep{SBB}.
Roughly 60\% of all wagons transported by \SBBCargo are generated by the \emph{single wagon load} transportation system (Wagenladungsverkehr or in short \WLV). Unlike unit train services where full trainloads are moved directly from origin to destination, WLV enables the transportation of individual wagons, making it ideal for customers with smaller freight volumes. These wagons follow a series of train segments across a network from receiving and formation stations to classification yards and back. These sequences are called \emph{\tcs} (see \fref{fig:transport_chain} for an example).
\SBBCargo's \WLV system comprises roughly 300 stations whereof five are used as major classification yards. Serving 10,000-15,000 wagons per week, \SBBCargo relies on stable and efficient algorithms for automatically generating \tcs within its IT system.

We briefly review past developments of \SBBCargo's transport chain computation.
Prior to 2017, production followed a leave-when-full and first-in-first-out methodology.
Train capacities were only considered at operational time, suspending wagons if necessary. We note that train capacities are used in Switzerland for weight and length to ensure that the train does not exceed the pulling capacities of its locomotives, can be over-taken by passenger trains and can obey the speed requirements prescribed by the infrastructure provider.
In 2017, \SBBCargo changed to a rigorous capacity management respecting train capacities at booking time. This was implemented by computing \tcs using a mapping with target zones: given a wagon at a current station and the target zone of its destination, this mapping would yield the subsequent station to travel to. To compute a full \tc, starting from the wagon's origin, the system iteratively returned the earliest train towards the subsequent station in the mapping still having sufficient capacity, until the destination was reached. Note that the mapping mechanism is closely related to the concept of so-called \emph{Leitwege}, also used by other major railway operators~\citep{Fügenschuh:2013}. Despite improvements, the approach remained rigid as \tcs were constrained by the predefined geographical mappings and fixed at booking time with limited flexibility to handle dynamic demands.

To modernize its IT and process landscape, \SBBCargo launched the project \emph{Greenfield} in 2019. Greenfield was intended to enable a more flexible production. In particular, the online nature of the bookings was mitigated using mathematical optimization of the \tcs.

For the management of the train network and the bookings, \SBBCargo decided to tailor the software Rail Cargo Management Solution by DXC Technology to its needs~\citep{Albrecht:2024, Albrecht:2024:2}. Its running instance is called \emph{ORCA} (ORder-to-CAsh). For the \tc computation, however, \SBBCargo decided to keep the development mostly in-house, implementing a software component called \emph{\TKModul} which communicates with ORCA to source all the required data for computing \tcs. The migration to the new systems was completed by the end of~2023.

The decision to keep the development of the \TKModul separated from ORCA had several advantages. First, we have direct control of the algorithms and their logic. This improved the understanding, level of detail, and applicability of the computation.
Second, we can also use the \TKModul for other purposes outside of the operational software, \eg, simulations of new production concepts and network designs.
Third, we integrated mathematical optimization models to react dynamically on fluctuations in demand. This is necessary because the train network is planned months in advance whereas the bulk of the ultimate transportation demands is not even known up to twelve hours in advance (cf. \fref{fig:bookings}). The algorithms can (re-)distribute wagons in real-time and break with the previous paradigm of simply catching the next train in a given direction.
We note that from theory, it is known that stochasticity in demands can degrade the performance of deterministic network designs and that dynamic capacity utilization and the ability to handle alternative routings is a prerequisite to reduce such effects~\citep{Wallace:2019}.

This paper describes the development and implementation of the \TKModul: we first describe the \tc search in \sref{tc search} including the required data structures. Then, we discuss the mathematical optimization and its online application in \sref{optimization}. \sref{adjustments} discusses additional requirements and operational details. In \sref{implementation}, we outline the technology stack, including deployment and interfaces. In \sref{operation}, we present some statistics about the operational system and conclude in \sref{conclusion}.

\section{Transport chain search}\label{tc search}

\subsection{Data structures}\label{sec:data_structures}

\begin{figure}[t]
	\centering
    \includegraphics{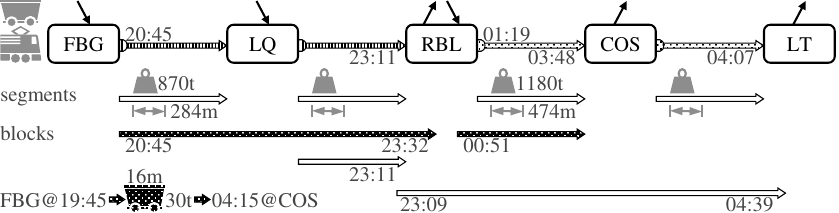}
	\caption{Example of train 50476 from Felsberg (FBG) via Landquart (LQ) to the shunting yard in Zurich (RBL), and train 50208 to Cossonay (COS) and shunting yard in Lausanne (LT). Wagons can be attached in FBG, LQ and RBL and detached in RBL, COS and LT. Depicted are selected train timings, capacities on the segments, boarding and deboarding times on the blocks. The filled blocks form a \tc for the depicted request from FBG to COS.}
	\label{fig:segments_and_blocks}
\end{figure}

We model a transportation demand as a \emph{request}, which represents one or several wagons or containers with specific weight and length, to be transported from its origin to its destination with earliest pickup and latest delivery times (see \fref{fig:segments_and_blocks}). Requests may also include attributes like customer identifiers, commodity codes (NHM codes), product types (\eg, time-critical), and maximum speed or coupling requirements.

The notion of capacities is ambiguous in train networks. We focus on a setting where the trains are already fixed (cf. \sref{optimization}) and have length and weight restrictions. In particular, we do not consider constraints on the number of trains in a timetable or along geographic corridors.

The train network is composed of trains with defined routes and intermediate stops. Each stop has specific purposes, such as allowing for the attachment or detachment of wagons, and trains have weight and length capacities between stops.
For transport chains, we use the concepts of \emph{blocks} and \emph{segments}. A segment represents the portion of a train route between two stops and manages capacities, while a block spans one or more segments and manages boarding and deboarding times (cf. \fref{fig:segments_and_blocks}).
The term block is closely related to the methodology of "blocking" used by North American railways to group wagons together for longer travel sections (see Chapter 13 in~\cite{Crainic:2021}), hence the naming.

For the simplest form of blocks, we generate one block between two stops along the same train that support attachment and detachment of wagons, respectively. Therefore, in \fref{fig:segments_and_blocks} we have two blocks for the train from FBG to RBL and two blocks for the train from RBL to LT. These blocks also account for key logistics processes, such as decoupling, inspection, and formation times. For example, the block timings at the shunting yards RBL and LT differ from their underlying train times because they include buffers for classification. Blocks can be restricted to specific requests through \emph{restrictions}, and transfers between blocks can be controlled via a \emph{transfer matrix}.

Finally, a \emph{\tc} is a sequence of blocks that are \emph{chainable}. This means they are geographically and timely consistent and all pairwise transfers between subsequent blocks are allowed. Additionally, chains must not form geographical loops. Hereby, the geography is solely checked on the origins and destinations of the blocks and not their contained segments. This means, physically, a wagon's \tc might very well loop in geography (\eg, imagine a first-mile pickup tour that travels forth and back along the same geographic path) but it cannot board or deboard at the same station twice.

A chain is considered \emph{valid} for a request if the pickup and delivery times are respected and the request satisfies all block restrictions and does not exceed segment capacities.

Transport chains are separated into \emph{required} and \emph{flexible} blocks. This concept is useful when a request is bound to a particular chain. The first part of required blocks represents facts that can no longer be changed, \eg, if the train has already departed. It can also be request-specific, \eg, if planners prescribe the \tc manually.
The second part of flexible blocks may still be changed and is subject to optimization.

Since the \WLV is not a self-contained system, there are requests and blocks outside of it. We refer to them as \emph{manual} and they are still communicated to the \TKModul since they may use the same resources. However, they are not included in chain searches and are fixed in the optimization.

\subsection{Search algorithms}\label{search algorithms}

For a given request, two routines are implemented in the \TKModul to search for valid \tcs. The first routine finds the single best valid chain, while the second enumerates all valid chains. Both routines can account for other requests' capacity usage or search for the request in isolation (as if there were no other in the network). Also, both routines respect required blocks and only complete chains in the flexible part.

To define the "best" chain, we use a tie-breaker which orders lexicographically by earliest departure, earliest arrival, fewest blocks, and earliest departures along the way.
The decision to give highest priority to the earliest departure was made before the implementation of restrictions and causes serious problems (cf. \sref{adjustments}).
Only recently, after the submission of this paper, we managed to change the tie-breaker to prioritize earliest arrival highest.
Note that the tie-breaker is solely used to improve the applicability of the solutions and not to improve the mathematical optimization models (\eg, by breaking symmetries).

We emphasize that just deciding if a valid chain exists for a request is NP-complete, because the blocks have a flexible neighborhood relation governed by a transfer matrix and valid chains must be free of geographical loops. To observe the complexity, think of the blocks as nodes in a graph and the destinations of the blocks as colors. Then, we are looking for a path through this graph (a sequence of chainable blocks) which does not repeat any color (geographically loop-free). This problem is called \emph{Rainbow s-t vertex connection} and is NP-complete by reduction from 3SAT~\citep{Chen:2011}. Note further that our graph on blocks is a directed acyclic graph because blocks move forward in time. But, the construction in \cite{Chen:2011} also uses a DAG and, hence, there is little hope in reducing the complexity.

Given the complexity, we employ an heuristic. Using breadth-first search (BFS) on blocks, we minimize the tie-breaker, starting from the request's origin at pickup time and ending at its destination at delivery time. If the found chain is not geographically loop-free, we discard it. In this case, we run the second routine, enumerating all chains, but stop once the first (best) chain candidate is found.

The second routine enumerates all valid chains using BFS as well. It operates on sequences of blocks - which are partial chains - in its priority queue. These can easily be kept loop-free. However, enumerating all chains has exponential runtime in the worst case, so we use limits to interrupt the search.

To avoid excessive computation times and memory usage, we impose limits on both algorithms. First, we limit the computation time, usually to five seconds. Second, we limit the search depth, usually to seven blocks per chain. Third, we limit the search frontier, usually to 10M partial chains. Finally, we limit the number of chains when enumerating all chains and stop once the best one hundred are found.

\section{Optimizing transport chain assignments}\label{optimization}

\subsection{Throughput optimization model}

\begin{figure}[t]
    \centering
    \includegraphics[width=0.73\linewidth]{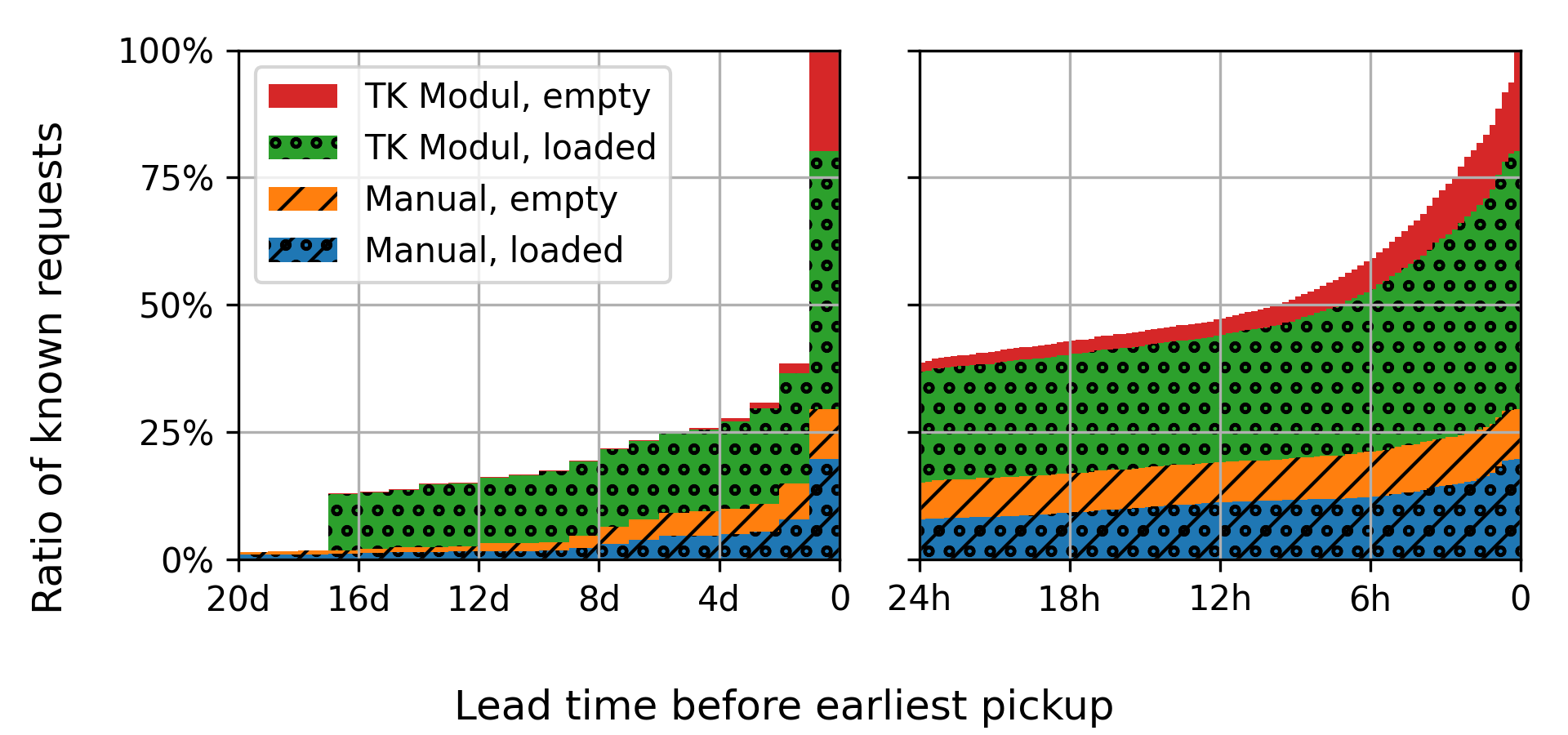}
    \caption{Lead time of the bookings before the requests' departures. Data from September 22nd to 28th, 2024. We omitted 9\% of the requests, for which there was no message prior to departure, as there is technically no lead time required and we cannot distinguish such showups from data inconsistencies.}
    \label{fig:bookings}
\end{figure}

Given that our planning horizon is much longer than the lead times (\fref{fig:bookings}), the production resources (\eg, locomotives, drivers, trains) are already fixed and cannot be changed. This also applies to the train capacities. The objective of the optimization is therefore to route as many requests as possible.

We use a path-based formulation below because \tcs are usually short, \ie, between three to six blocks long. Hence, one may enumerate chains directly, rather than having to resort to column generation techniques or using flow-based formulations.

Let $R$ be the set of requests, $C_r$ be the set of chains that are valid for request~$r$, $B_c$ the blocks along chain~$c$, and $S_b$ the segments along block~$b$.
We simplify notation by using $B$ as the set of all blocks ($\bigcup_{R\in R}\bigcup_{c\in C_r} B_c$) and $S$ the set of all segments ($\bigcup_{b\in B} S_b$).
To further ease notation we will denote the back-pointers by $R_c$ for all requests that use chain~$c$, $C_b$ for all chains that contain block~$b$ and $B_s$ for all blocks that contain segment~$s$. That means, for instance, that $B_s = \{b\in B\ |\ s \in S_b\}$.

Each request has capacity requirements $\ubar{c}_r^i$ for $i\in\{\textnormal{weight}, \textnormal{length}\}$. Similarly, each segment has maximum capacities $\overline{C}_s^i$.
We let $x_{r,c}$ denote the indicator if request~$r$ is routed over its chain $c\in C_r$.
Now, clearly we would like to maximize (a weighted) sum of the $x_{r,c}$ under the constraints that all capacities are satisfied, \ie,
$\sum_{b\in B_s}\sum_{c\in C_b}\sum_{r\in R_c} \ubar{c}_r^i \cdot x_{r,c} \le \overline{C}_s^i$.

But, life is not that simple. 
In practice, we are also given capacity reservations on the blocks. These are used to reserve space for anticipated requests whose bookings are late, usually combined with some restrictions on the blocks to tailor the space reservation to a precise target group of requests.
To accommodate such capacity reservations we denote them by $\text{res}_b^i$ and add variable capacities $\text{cap}_b^i$ for each block~$b$.

Furthermore, some of the requests in $R$ might already have an assigned chain. We need to ensure that these requests, denoted by $R_A$, will keep a routing even though they may be rerouted on different chains. We emphasize that the pickup and delivery times are hard constraints enforced in the \tc search. That means even if a request is rerouted to a new \tc it still obeys its time limits.

Then, the optimization model is
\begin{align}
    \text{\customlabel{opt model}{(OPT)}}\quad\max && \sum_{r\in R} \sum_{c\in C_r} p_{r,c} \cdot x_{r,c} \nonumber\\
    \text{s.t.} && \sum_{c\in C_r} x_{r,c} & \le 1 && \forall r\in R\setminus R_A \label{eq:routing}\\
    && \sum_{c\in C_r} x_{r,c} & = 1 && \forall r\in R_A \label{eq:fixed}\\
    && u_b^i + \sum_{c\in C_b}\sum_{r\in R_c} \ubar{c}_r^i \cdot x_{r,c} & \le \text{cap}_b^i && \forall i, b\in B \label{eq:block capacity}\\
    && \sum_{b\in B_s} \text{cap}_b^i &\le \overline{C}_s^i && \forall i, s\in S \label{eq:segment capacity}\\
    && \text{res}_b^i &\le \text{cap}_b^i && \forall i, b\in B \label{eq:reservations}\\
    && x_{r,c} & \in \{0,1\} && \forall r\in R, c\in C_r \nonumber \\
    && \text{cap}_b^i & \in \mathbb{Z}_{\ge 0} && \forall i, b\in B \nonumber
\end{align}
Here, $p_{r,c}$ is a coefficient that is used to implement priorities for routing certain requests as well as the preference for chains that have a low objective in the tie-breaker (as defined in \sref{search algorithms}).
Constraints~\eqref{eq:routing} and~\eqref{eq:fixed} ensure that the requests obtain at most one chain and the requests in $R_A$ are certainly (re)routed.
Constraint~\eqref{eq:block capacity} ensures that the (variable) capacity of a block is at least its utilization by the requests. Thereby,~$u_b^i$ is the capacity utilization of potential requests that should not be touched, \ie, are outside of $R$. 
Constraint~\eqref{eq:segment capacity} ensures consistency with the segment capacities and~\eqref{eq:reservations} ensures that each block takes up its reserved space.

The attentive reader might argue, that this model may be infeasible if the reservations~$\text{res}_b^i$ are set too high in \eqref{eq:reservations} for the capacity requirement of~$\overline{C}_s^i$ in \eqref{eq:segment capacity}. The \TKModul should neither introduce nor remove over-booking, thus the capacity of a segment is defined to be the maximum of its planned capacity and the sum of all requests and reservations using it.

In production our models consist mostly of only a dozen variables and constraints, as the model is used with a single request 94\% of the time. The largest model we saw in production had 4,469 variables and 1,431 constraints for assigning 117 requests, and was still built and solved to optimality in 1.3s. In simulations we solved models for 367K requests in 5m to sufficient precision. Given solving was never a performance issue, we have not bench-marked the model or optimized the solver settings.

\subsection{Application in online mode}\label{online-application}

We have not yet specified where the set of requests $R$ comes from.
In \SBBCargo's productive system, requests are not known in advance but the bookings arrive in an online fashion. \fref{fig:bookings} shows the lead times, \ie, the time difference between the booking and the departure of a request. There is a notable cluster around two weeks due to 'contingents' that reserve capacity for anticipated requests. However, less than half of the requests are known more than 12 hours in advance, requiring the \TKModul to make online and real-time decisions about accepting and routing requests. Once a request has been accepted, this decision cannot be revised. However, requests can be rerouted to free up space for incoming request bookings, as long as they are still delivered in time.

A simple approach is to route requests sequentially, using a greedy strategy: if there is enough capacity, route the request; if not, reject it. However, this leads to congestion and suboptimal results. Early experiments in an offline setting, where the order of the requests can be adjusted, showed that a greedy approach missed 7-12\% of requests compared to an optimal solution, especially when train capacities were tight. Performance varied based on how requests were ordered, with random ordering - representing the online setting - producing the worst outcomes and the highest inconsistency in routing results.

While in theory, a full optimization over all requests and \tcs, could be run at each request's arrival, this would be computationally expensive and result in frequent changes of \tcs. Therefore, we opted for a local optimization approach which is computationally scalable, maintains routing stability, and is easier to debug.

\begin{varalgorithm}{Assignment}
    \caption{Pseudo-code for algorithmic optimization of \tc assignments.}
    \label{alg:online assignment}
    \begin{algorithmic}
        \Require Set of requests $R$ %, chains $C$, blocks $B$ and segments $S$ including their current capacity utilization
        \Function{Try Optimization}{$R$, $C$} \Comment{shortcut}
            \State run \ref{opt model} for $R$ with the $C_r$
            \State let $R_1$ be the routed and $R_2$ the non-routed requests
            \If{$R_1=R$ ($\Leftrightarrow R_2 = \emptyset$)}
                return (globally) the obtained assignment
            \EndIf
        \EndFunction
        \State
        \State $C_r \gets \emptyset \quad \forall r\in R$ \Comment{set of chain candidates for each request}
        \For{$r\in R$} \Comment{step 1: find best chains in isolation}
            \State find (single) best chain for $r$ in isolation and add it to $C_r$
            \If{none exists}
                label $r$ as not routable and remove it from $R$
            \EndIf
        \EndFor
        \State \Call{Try Optimization}{$R$, $C$}
        \State $\text{neighborhood\_search} \gets \text{True}$
        \For{$r\in R_2$} \Comment{step 2: add chain respecting capacities}
            \State find (single) best chain for $r$ respecting capacities and add it to $C_r$
            \State if none exists, neighborhood\_search $\gets$ False and break \Comment{move to step 4}
        \EndFor
        \If{neighborhood\_search}
            \Call{Try Optimization}{$R$, $C$}
        \EndIf
        \For{$r\in R$} \Comment{step 3: find all chains}
            \State find all valid chains for $r$ respecting capacities and add them to $C_r$
        \EndFor
        \If{neighborhood\_search}
            \Call{Try Optimization}{$R$, $C$}
        \EndIf
        \State find the current neighborhood $\overline{R}$ w.r.t. $C$ \Comment{step 4: optimize including the neighborhood}
        \For{$r\in \overline{R}\setminus{R}$}
            \State find all valid chains for $r$ respecting capacities and set $C_r$ to be this set
        \EndFor
        \State \Call{Try Optimization}{$R$, $C$}
        \State return found assignment for requests in $R_1$
    \end{algorithmic}
\end{varalgorithm}

The algorithm for online chain assignment is outlined in \aref{alg:online assignment}. It uses a function $\Call{Try Optimization}{}$ which runs the mathematical optimization \ref{opt model} for a specified set of the requests and their respective chain candidates. Note that the first step, where we search for each request in isolation, is unnecessary for maximizing throughput and could be skipped. However, in a productive system it is important to distinguish if a valid \tc exists at all or if there is no remaining capacity.

In step 4, we re-optimize \tcs over a neighborhood. The \emph{neighborhood} of a set of chains~$C$ are all requests whose current chain shares a segment with any of the chains in~$C$.
While a broader optimization could be run across all requests, this would increase runtime and result in less stable routing decisions. Note, that we could also use a second, third, or deeper neighborhoods iteratively in step~4. But for operational tractability, we usually limit the neighborhoods to the first one. The median number of requests in the first neighborhood was 63, while rare cases with several hundred requests in the first neighborhood ran into time limits when finding all \tcs.

\begin{figure}[t]
    \centering
    \includegraphics[width=\textwidth]{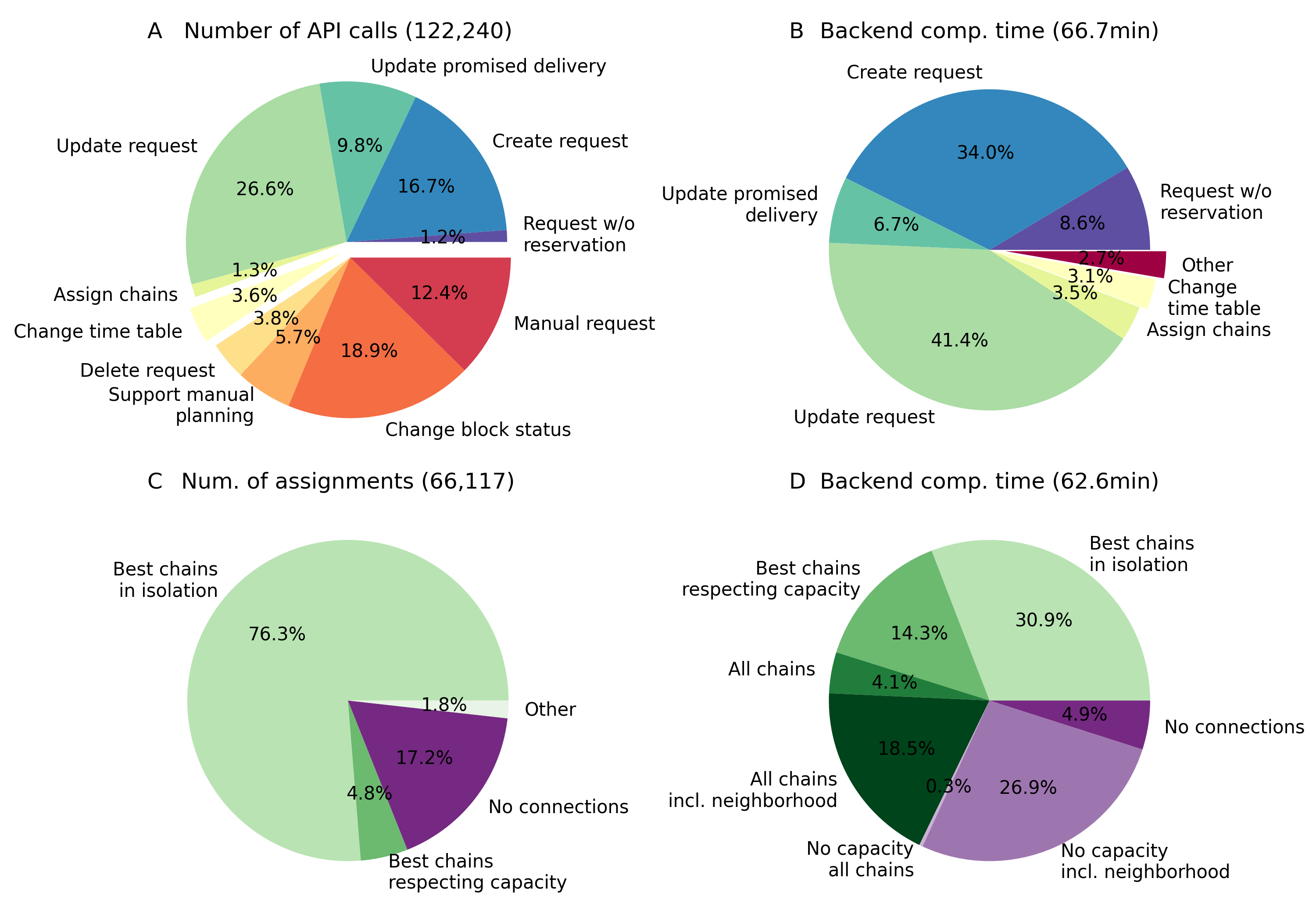}
    \caption{Data from September 22nd to 28th, 2024. A) Clustered calls to the \TKModul's API. The upper half (purple to green) also triggered \aref{alg:online assignment}. B) Computational time in the backend by cluster. C) The outcomes of the assignments. D) Computational time in the backend used by the assignments by outcome.}
    \label{fig:operation}
\end{figure}

To ensure latency, \aref{alg:online assignment} can be stopped at each intermediate step if time runs out (yielding the different outcomes in \fref{fig:operation}). Further, chains are cached to reuse them for similar requests.

\section{Additional details for operational application}\label{adjustments}

\subsection{Customizing the search for productional feasibility}

To obtain practical \tcs, it is important to account for transportation restrictions. Therefore, blocks can include restrictions that filter which requests can use them. Each restriction targets a specific request attribute (\eg, origin, receiver, NHM codes) detailing allowed or forbidden values.

The construction of \tcs from blocks uses transfers between blocks. By default, a transfer is feasible if the respective boarding and deboarding times align. More fine grained configuration is possible with a transfer matrix using \emph{connections}. Available types are forbidden connections, extra connections and exclusive connections. Extra connections can be quick transfers between trains that bypass the hump in a shunting yard. Exclusive connections can, \eg, be used to feed all the block's requests to a subsequent one.

During operations at a shunting yard, wagons are often set aside for later train formation phases. This is done by pushing them over the hump into a collection track whose contents are later pushed a second time over the hump. To model these operations, we introduced phase connectors as additional blocks.

Blocks have a wide modeling power. For instance, they are used to help enforce formation groups which are groups within the train sharing the same destination. They could also span the segments of multiple trains.
Furthermore, the \tcs usually start at the first and end at the last train. One may also want to include the first and last mile, \ie, the pickups and deliveries (usually by diesel locomotive) at the customer service point. These can also be modeled as blocks, thereby enabling additional logic between service trips and first/last trains via the transfer matrix. 
Unfortunately, ORCA only generates blocks for a single train (the simple form in \sref{sec:data_structures}), but we use these possibilities in tactical simulations.

\subsection{Technical adjustments to mappings and workflow}

To chain two blocks, the destination of the former must match the origin of the latter. In practice, however, a shunting yard might comprise several operation points, \eg, arrival, classification and departure groups. To span the distance between these technically different locations, we use groups of operation points within which transfers are possible. Similarly, requests live on a commercial layer whereas segments and blocks live on an operational layer. The gap is bridged by a mapping from commercial to operational stations.

When the train network (including reservations, restrictions, and connections) changes, the affected chains are checked for feasibility. This is done block by block and they are truncated after the last consecutive feasible block. New chain searches only try to complete in the flexible part of the \tc. If this is not successful, the partial chains are still maintained to keep the requests moving in the right direction.

When booking, the customer can choose a service time at his location. Subsequently, ORCA looks up the next service window and uses its end as pickup for the resulting request. The delivery time for the request is deduced by the product, \eg, express transports obtain a maximal time window of 29 hours. Then, the request is forwarded to \TKModul, which searches and returns a \tc, and the arrival time of the chain is communicated back to the customer. Note that the punctuality of the request - a major KPI at the Swiss national railways - is measured against this earliest estimate.
Therefore, to enforce the delivery time, ORCA sends an update for the same request shortening the time window to the promised delivery time immediately after obtaining the first \tc from the \TKModul. Unfortunately, this workflow has serious disadvantages. First, the promised delivery times significantly reduce the number of valid chain candidates up to the point of supporting only one chain. Therefore, from an algorithmic perspective, the optimization degrades to a greedy assignment (cf. \sref{online-application}). %And to mitigate congestion, practitioners introduce and maintain restrictions to strongly guide the search to avoid bottlenecks.
Second, the \tc computation can yield operationally infeasible results mainly due to the previous tie breaker which prioritized earliest departure. This may cause requests to travel long ways because they follow the first train leaving at their current location. Of course, this "Tour de Suisse"-phenomenon is kept within limits by tuning the configurations.
Third, the IT systems of the shunting yards are not directly coupled with ORCA. Therefore, they may override the planned \tcs at their own discretion. This helps if the computed \tc is operationally infeasible but it produces deviations from the promised delivery.

\subsection{Over-steering the \TKModul}

In \sref{sec:data_structures}, we introduced the concept of required blocks to represent parts of a \tc that can no longer be changed. Such required blocks can also be used for manual planning. Manual chains are required, for instance, to make sure equipment arrives at construction sites facing the correct direction. Geographic gaps, loops or infeasible transitions within required blocks are accepted by the \TKModul.

Both, operational deviations from planned capacities and manual planning can lead to over-booking. The \TKModul accepts over-booking even on bookable blocks and does not correct it when re-assigning chains. However, when the capacities in the planned train network are reduced, the \TKModul will search for new chains and avoid over-booking on the flexible blocks. 

\section{Implementation, technology stack, and development}\label{implementation}

In a productive setting, \tcs must be computed in real-time, which raises the desire for parallelization.
However, segment capacities are global properties that cannot be exceeded and, thus, conflicts between requests - potentially processed in parallel - must be avoided. Unfortunately, for parallelization, a time-wise separation is impossible since lead times are small (cf. \fref{fig:bookings}) and, therefore, the travel time windows of most requests overlap. Also, a geographical separation is impractical because the requests' origins and destinations are everywhere in Switzerland and their \tcs meet up in only a few shunting yards. Therefore, \tc computations run sequentially in a single process.

To ensure response times below five seconds, the \TKModul must have fast access to all business objects described in \sref{tc search} and \sref{adjustments}. This led to the decision to store all these objects in-memory in a so-called \emph{state}. To give an example of its size, the state in the production environment on October 26, 2024, contained 61k segments, 75k blocks, 42k requests, and 3k connections.

The in-memory design has disadvantages because it necessitates a complex synchronization mechanism between ORCA, where the train network and the requests are managed, and \TKModul, where the state resides.
To mitigate complexity, we implemented the \TKModul as a server without persistent storage that, bi-directionally, communicates with ORCA through HTTP. When starting up, it receives an initial message with all relevant business objects and master data and subsequently obtains frequent updates, during peak hours with an average of 1,200 messages every five minutes.

Planners at \SBBCargo need to compute \tcs just for their information and they need to test the feasibility of manually planned \tcs. To support these tasks, we implemented endpoints that dry run the computation, without reserving capacities, and also collect reasons for infeasibility, \eg, non-matching restrictions, insufficient capacities, or infeasible transfers between the blocks.

Changes in the train network may require multiple messages, each affecting the current routing, which in turn affect the \tcs and may require new \tc computations. To avoid these immediate triggers, most endpoints provide a flag to \emph{not} yet run the \tc computation.
It is eventually triggered, once the update on the train network and its batch of messages has been processed.

To keep the API responsive, the \tc searches and optimization run (sequentially) in a backend queue. This keeps the endpoints open to receive further messages. The backend queue also implements priorities, which favor operations triggered by the end users. For instance, if large parts of the train network are modified, this may impact lots of requests whose new chain searches are then carried out on the side.

Technologically, the \TKModul is implemented as a Python library. It uses Pydantic~\citep{Pydantic} for validation of the business objects and messages on the API to ORCA. It has an executable HTTP server, using FastAPI~\citep{FastAPI} which runs in its own Docker container without persistent storage, hosted on an Openshift cluster of the SBB. Logs are collected in Splunk, where we built a detailed dashboard for monitoring. We use SCIP~\citep{SCIP} as mixed-integer programming solver.

To date, four mathematicians, one physicist, and one industrial engineer have made significant contributions to the development of the \TKModul. While the frequent rotation presented challenges, it also had benefits: it drove continuous improvement of the most complex parts of the code, ensured experienced reviewers were always available, and necessitated thorough testing.

Integration testing with both \TKModul and ORCA was notoriously difficult. Issues surfaced almost exclusively in production, as the data in the test and integration environments did not cover the full complexity. Performance tests proved particularly unreliable, leaving uncertainty about whether our architecture would be able to handle the rapid increase in load during the steep migration ramp-up \citep{Albrecht:2024}.

\section{Operation and performance}\label{operation}

To illustrate how the \TKModul operates, we analyzed the logs from the same week as in \fref{fig:bookings}, which is representative for a week without larger changes to the train network. During that week, the \TKModul's API was called 122,240 times, which activated the backend queue for 67 minutes. About 56\% of these calls triggered \aref{alg:online assignment}, using about 94\% of the time the backend queue was active. Of the 66,117 assignments more than 98\% only required the cheaper breadth-first search, however, the remaining 1.8\% required about half the computation time in the backend.

We further logged the computational time spent in the most time consuming functions. We found the backend to spend 30min to enumerate all chains, 25min in the breadth-first search, and 5min 17s in the optimization, of which only 47s are spent solving the almost 60k models, the rest building them. Another 9min was used by the synchronous endpoints to support planning, one minute for the breadth-first search and 8min to enumerate all chains. To deepen our understanding about where the computation time is spent, we simulated 24h of messages to the \TKModul locally with a profiler and observed that the majority of the time is spent listing the next possible blocks, \eg, collecting visited stations to avoid geographic loops and looking up exclusive connections, functionality we built with simplicity and not with performance in mind.

\section{Conclusions, lessons learned, and future work}\label{conclusion}

One year after the successful migration to the new systems, we can draw a positive conclusion. The migration proceeded with the \TKModul rarely being the source of issues and we were able to overcome the usual implementation challenges of mathematical optimization in railways (see \cite{Liebchen:2019}). Further, operational costs are drastically lower, as the old systems ran on an IBM mainframe. However, comparing reliability between the two systems remains difficult, as they measure punctuality differently and A/B testing is impossible in production.

Keeping the \TKModul in-house turned out to be a good decision. The \tc computation is at the heart of the \WLV system and \SBBCargo thereby ensured that it keeps in-depth knowledge about this computation. Further, it required ORCA to be open with data, revealing bugs within ORCA and its APIs early. This also facilitated migration of data between test and production environments.

Using the \TKModul for other purposes is easily possible because it is implemented as a Python library where the core algorithms are separated from the API to ORCA. \SBBCargo now also uses it for various types of simulations, comparing \tcs across different train networks, evaluating tactical network designs, or assessing the impact of changes to the objective function.

Unfortunately, \SBBCargo suffers from a decade-long decrease in revenues that causes large monetary deficits. The steady decrease in the number of wagons compared to the stable train capacities reduces congestion. Therefore, the bulk of the \tc computation can be handled by simple searches without the need for mathematical optimization. But, this mathematical component is also used outside of the operational software. After optimizing the train network for resources on a tactical level (by different models), the \TKModul is used to simulate effective capacity usage.

Moving to a fully dynamic capacity management proved to be too complex. In part this is due to the complexity in railway production where planners traditionally envision a single possible \tc per request. Handling alternatives is often deemed infeasible. The new system forces planners to configure a set of rules rather than plan direct routes, which requires change management.
Further, leveraging the full optimization potential would have required a major overhaul of the booking process. To date, deliveries are promised at booking time and used to measure punctuality. Delaying this information slightly, as done by most parcel services, would require adjustments in several IT systems and processes.

There is still ample room for improvement. The outcomes from the well established previous system were more predictable. As transport chains followed strict geographic routes, wagons would simply catch the next train heading in the right direction. In contrast, the flexibility of the new system can sometimes lead to infeasible detours, requiring manual intervention.
To this end, we recently changed the tie-breaker in the prioritization of a request's chains to remove earliest departure as the primary objective. While a tie-breaker should not matter in theory, it is rather delicate in practice: We almost exclusively see the "best" \tc in production and thus the \TKModul's configuration is much less tested beyond.
Another future improvement might be including the block-to-train problem (see \cite{Harrod:2011} for a review) to dynamically minimize shunting movements. We currently have to hard-wire some connections in the transfer matrix to ensure feasibility for the shunting teams. However, this would also require the ability to adopt short-term changes at the shunting yards.
Finally, the TK Modul would be an appropriate system to automate the distribution of empty wagons, which we currently do manually in production and crudely approximate in tractical simulations.

\bibliographystyle{apacite}
\urlstyle{rm}
\bibliography{tkmodul}

\begin{thebibliography}{}

\bibitem [\protect \citeauthoryear {%
Albrecht%
, Skujat%
, L\"{o}rincze%
\BCBL {}\ \BBA {} vom Hagen%
}{%
Albrecht%
\ \protect \BOthers {.}}{%
{\protect \APACyear {2024}}%
}]{%
Albrecht:2024}
\APACinsertmetastar {%
Albrecht:2024}%
\begin{APACrefauthors}%
Albrecht, T.%
, Skujat, N.%
, L\"{o}rincze, G.%
\BCBL {}\ \BBA {} vom Hagen, T.%
\end{APACrefauthors}%
\unskip\
\newblock
\APACrefYearMonthDay{2024}{September}{}.
\newblock
{\BBOQ}\APACrefatitle {{Agile Modernisierung der Order-to-Cash-Softwarelandschaft bei SBB Cargo}} {{Agile Modernisierung der Order-to-Cash-Softwarelandschaft bei SBB Cargo}}.{\BBCQ}
\newblock
\APACjournalVolNumPages{Eisenbahntechnische Rundschau}{}{}{}.
\PrintBackRefs{\CurrentBib}

\bibitem [\protect \citeauthoryear {%
Albrecht%
, Skujat%
, Lörincze%
\BCBL {}\ \BBA {} vom Hagen%
}{%
Albrecht%
\ \protect \BOthers {.}}{%
{\protect \APACyear {2025}}%
}]{%
Albrecht:2024:2}
\APACinsertmetastar {%
Albrecht:2024:2}%
\begin{APACrefauthors}%
Albrecht, T.%
, Skujat, N.%
, Lörincze, G.%
\BCBL {}\ \BBA {} vom Hagen, T.%
\end{APACrefauthors}%
\unskip\
\newblock
\APACrefYearMonthDay{2025}{}{}.
\newblock
{\BBOQ}\APACrefatitle {A planning framework for high automation of rail cargo order processing} {A planning framework for high automation of rail cargo order processing}.{\BBCQ}
\newblock
\APACjournalVolNumPages{11th International Conference on Railway Operations Modelling and Analysis RailDresden}{}{}{}.
\PrintBackRefs{\CurrentBib}

\bibitem [\protect \citeauthoryear {%
Bestuzheva%
\ \protect \BOthers {.}}{%
Bestuzheva%
\ \protect \BOthers {.}}{%
{\protect \APACyear {2023}}%
}]{%
SCIP}
\APACinsertmetastar {%
SCIP}%
\begin{APACrefauthors}%
Bestuzheva, K.%
, Besan\c{c}on, M.%
, Chen, W\BHBI K.%
, Chmiela, A.%
, Donkiewicz, T.%
, van Doornmalen, J.%
\BDBL {}Witzig, J.%
\end{APACrefauthors}%
\unskip\
\newblock
\APACrefYearMonthDay{2023}{June}{}.
\newblock
{\BBOQ}\APACrefatitle {{Enabling Research through the SCIP Optimization Suite 8.0}} {{Enabling Research through the SCIP Optimization Suite 8.0}}.{\BBCQ}
\newblock
\APACjournalVolNumPages{ACM Transactions on Mathematical Software}{49}{2}{}.
\newblock
\begin{APACrefDOI} \doi{10.1145/3585516} \end{APACrefDOI}
\PrintBackRefs{\CurrentBib}

\bibitem [\protect \citeauthoryear {%
Chen%
, Li%
\BCBL {}\ \BBA {} Shi%
}{%
Chen%
\ \protect \BOthers {.}}{%
{\protect \APACyear {2011}}%
}]{%
Chen:2011}
\APACinsertmetastar {%
Chen:2011}%
\begin{APACrefauthors}%
Chen, L.%
, Li, X.%
\BCBL {}\ \BBA {} Shi, Y.%
\end{APACrefauthors}%
\unskip\
\newblock
\APACrefYearMonthDay{2011}{}{}.
\newblock
{\BBOQ}\APACrefatitle {The complexity of determining the rainbow vertex-connection of a graph} {The complexity of determining the rainbow vertex-connection of a graph}.{\BBCQ}
\newblock
\APACjournalVolNumPages{Theoretical Computer Science}{412}{35}{4531-4535}.
\newblock
\begin{APACrefDOI} \doi{10.1016/j.tcs.2011.04.032} \end{APACrefDOI}
\PrintBackRefs{\CurrentBib}

\bibitem [\protect \citeauthoryear {%
Colvin%
\ \protect \BOthers {.}}{%
Colvin%
\ \protect \BOthers {.}}{%
{\protect \APACyear {2024}}%
}]{%
Pydantic}
\APACinsertmetastar {%
Pydantic}%
\begin{APACrefauthors}%
Colvin, S.%
, Jolibois, E.%
, Ramezani, H.%
, Garcia~Badaracco, A.%
, Dorsey, T.%
, Montague, D.%
\BDBL {}Hall, A.%
\end{APACrefauthors}%
\unskip\
\newblock
\APACrefYearMonthDay{2024}{9}{}.
\newblock
\APACrefbtitle {Pydantic.} {Pydantic.}
\newblock
\begin{APACrefURL} \url{https://docs.pydantic.dev/latest/} \end{APACrefURL}
\PrintBackRefs{\CurrentBib}

\bibitem [\protect \citeauthoryear {%
Crainic%
, Gendreau%
\BCBL {}\ \BBA {} Gendron%
}{%
Crainic%
\ \protect \BOthers {.}}{%
{\protect \APACyear {2021}}%
}]{%
Crainic:2021}
\APACinsertmetastar {%
Crainic:2021}%
\begin{APACrefauthors}%
Crainic, T\BPBI G.%
, Gendreau, M.%
\BCBL {}\ \BBA {} Gendron, B.%
\end{APACrefauthors}%
\unskip\
\newblock
\APACrefYear{2021}.
\newblock
\APACrefbtitle {Network Design with Applications to Transportation and Logistics} {Network design with applications to transportation and logistics}.
\newblock
\APACaddressPublisher{}{Springer International Publishing}.
\newblock
\begin{APACrefDOI} \doi{10.1007/978-3-030-64018-7} \end{APACrefDOI}
\PrintBackRefs{\CurrentBib}

\bibitem [\protect \citeauthoryear {%
{Federal Statistical Office}%
}{%
{Federal Statistical Office}%
}{%
{\protect \APACyear {2023}}%
}]{%
FSO}
\APACinsertmetastar {%
FSO}%
\begin{APACrefauthors}%
{Federal Statistical Office}.%
\end{APACrefauthors}%
\unskip\
\newblock
\APACrefYearMonthDay{2023}{}{}.
\newblock
\APACrefbtitle {Statistics, Mobility and transport, Goods transport.} {Statistics, mobility and transport, goods transport.}
\newblock
\begin{APACrefURL} \url{https://www.bfs.admin.ch/bfs/en/home/statistics/mobility-transport/goods-transport.html} \end{APACrefURL}
\PrintBackRefs{\CurrentBib}

\bibitem [\protect \citeauthoryear {%
F\"{u}genschuh%
, Homfeld%
\BCBL {}\ \BBA {} Sch\"{u}lldorf%
}{%
F\"{u}genschuh%
\ \protect \BOthers {.}}{%
{\protect \APACyear {2015}}%
}]{%
Fügenschuh:2013}
\APACinsertmetastar {%
Fügenschuh:2013}%
\begin{APACrefauthors}%
F\"{u}genschuh, A.%
, Homfeld, H.%
\BCBL {}\ \BBA {} Sch\"{u}lldorf, H.%
\end{APACrefauthors}%
\unskip\
\newblock
\APACrefYearMonthDay{2015}{}{}.
\newblock
{\BBOQ}\APACrefatitle {Single-Car Routing in Rail Freight Transport} {Single-car routing in rail freight transport}.{\BBCQ}
\newblock
\APACjournalVolNumPages{Transportation Science}{49}{1}{130-148}.
\newblock
\begin{APACrefDOI} \doi{10.1287/trsc.2013.0486} \end{APACrefDOI}
\PrintBackRefs{\CurrentBib}

\bibitem [\protect \citeauthoryear {%
Harrod%
\ \BBA {} Gorman%
}{%
Harrod%
\ \BBA {} Gorman%
}{%
{\protect \APACyear {2011}}%
}]{%
Harrod:2011}
\APACinsertmetastar {%
Harrod:2011}%
\begin{APACrefauthors}%
Harrod, S.%
\BCBT {}\ \BBA {} Gorman, M\BPBI F.%
\end{APACrefauthors}%
\unskip\
\newblock
\APACrefYearMonthDay{2011}{}{}.
\newblock
{\BBOQ}\APACrefatitle {Operations Research for Freight Train Routing and Scheduling} {Operations research for freight train routing and scheduling}.{\BBCQ}
\newblock
\APACjournalVolNumPages{Wiley Encyclopedia of Operations Research and Management Science}{}{}{}.
\newblock
\begin{APACrefDOI} \doi{10.1002/9780470400531.eorms1014} \end{APACrefDOI}
\PrintBackRefs{\CurrentBib}

\bibitem [\protect \citeauthoryear {%
Liebchen%
\ \BBA {} Sch\"{u}lldorf%
}{%
Liebchen%
\ \BBA {} Sch\"{u}lldorf%
}{%
{\protect \APACyear {2019}}%
}]{%
Liebchen:2019}
\APACinsertmetastar {%
Liebchen:2019}%
\begin{APACrefauthors}%
Liebchen, C.%
\BCBT {}\ \BBA {} Sch\"{u}lldorf, H.%
\end{APACrefauthors}%
\unskip\
\newblock
\APACrefYearMonthDay{2019}{}{}.
\newblock
{\BBOQ}\APACrefatitle {A collection of aspects why optimization projects for railway companies could risk not to succeed – A multi-perspective approach} {A collection of aspects why optimization projects for railway companies could risk not to succeed – a multi-perspective approach}.{\BBCQ}
\newblock
\APACjournalVolNumPages{Journal of Rail Transport Planning \& Management}{11}{}{100149}.
\newblock
\begin{APACrefDOI} \doi{10.1016/j.jrtpm.2019.100149} \end{APACrefDOI}
\PrintBackRefs{\CurrentBib}

\bibitem [\protect \citeauthoryear {%
Ramírez%
}{%
Ramírez%
}{%
{\protect \APACyear {2024}}%
}]{%
FastAPI}
\APACinsertmetastar {%
FastAPI}%
\begin{APACrefauthors}%
Ramírez, S.%
\end{APACrefauthors}%
\unskip\
\newblock
\APACrefYearMonthDay{2024}{9}{}.
\newblock
\APACrefbtitle {FastAPI.} {Fastapi.}
\newblock
\begin{APACrefURL} \url{https://fastapi.tiangolo.com} \end{APACrefURL}
\PrintBackRefs{\CurrentBib}

\bibitem [\protect \citeauthoryear {%
{SBB Facts and Figures}%
}{%
{SBB Facts and Figures}%
}{%
{\protect \APACyear {2023}}%
}]{%
SBB}
\APACinsertmetastar {%
SBB}%
\begin{APACrefauthors}%
{SBB Facts and Figures}.%
\end{APACrefauthors}%
\unskip\
\newblock
\APACrefYearMonthDay{2023}{}{}.
\newblock
\APACrefbtitle {Freight Services.} {Freight services.}
\newblock
\begin{APACrefURL} \url{https://reporting.sbb.ch/} \end{APACrefURL}
\PrintBackRefs{\CurrentBib}

\bibitem [\protect \citeauthoryear {%
Wang%
, Crainic%
\BCBL {}\ \BBA {} Wallace%
}{%
Wang%
\ \protect \BOthers {.}}{%
{\protect \APACyear {2019}}%
}]{%
Wallace:2019}
\APACinsertmetastar {%
Wallace:2019}%
\begin{APACrefauthors}%
Wang, X.%
, Crainic, T\BPBI G.%
\BCBL {}\ \BBA {} Wallace, S\BPBI W.%
\end{APACrefauthors}%
\unskip\
\newblock
\APACrefYearMonthDay{2019}{}{}.
\newblock
{\BBOQ}\APACrefatitle {Stochastic Network Design for Planning Scheduled Transportation Services: The Value of Deterministic Solutions} {Stochastic network design for planning scheduled transportation services: The value of deterministic solutions}.{\BBCQ}
\newblock
\APACjournalVolNumPages{INFORMS Journal on Computing}{31}{1}{153-170}.
\newblock
\begin{APACrefDOI} \doi{10.1287/ijoc.2018.0819} \end{APACrefDOI}
\PrintBackRefs{\CurrentBib}

\end{thebibliography}

\end{document}